\documentclass[11pt,english]{amsart}
\usepackage[T1]{fontenc}
\usepackage[latin9]{inputenc}
\usepackage{amsbsy}
\usepackage{amstext}
\usepackage{amsthm}

\makeatletter
\usepackage{amsmath,amsfonts,amssymb,amsthm,epsfig}

\usepackage{color}
\usepackage[final,allcolors=blue,colorlinks=true]{hyperref}

\def\R{\mathbb R}

\newcommand{\dist}{\text{\rm dist}}
\newcommand{\medint}{-\kern -,375cm\int}         
\newcommand{\medintinrigo}{-\kern -,315cm\int}

\newcommand{\dive}{{\rm div}}

\numberwithin{equation}{section}
\newtheorem{theorem}{Theorem}[section]
\newtheorem*{theorem*}{Theorem}  

\newtheorem*{conclusion*}{Conclusin}

\newtheorem*{corollary*}{Corollary}

\newtheorem*{lemma*}{Lemma}

\newtheorem*{notation*}{Notation}
\newtheorem{problem}[theorem]{Problem}

\newtheorem*{proposition*}{Proposition}
\newtheorem{remark}[theorem]{Remark}
\newtheorem*{remark*}{Remark}

\newtheorem*{example*}{Example}                
\theoremstyle{definition}

\makeatother

\usepackage{babel}
\begin{document}
\title[]{The conservation law approach in geometric PDEs}

   \author[C.-Y. Guo and C.-L. Xiang]{Chang-Yu Guo and Chang-Lin Xiang$^\ast$}

\address[Chang-Yu Guo]{Research Center for Mathematics and Interdisciplinary Sciences, Shandong University, Qingdao 266237, People's Republic of China and Frontiers Science Center for Nonlinear Expectations, Ministry of Education, Qingdao, People's Republic of China}
\email{changyu.guo@sdu.edu.cn}

\address[Chang-Lin Xiang]{Three Gorges Mathematical Research Center, China Three Gorges University,  443002, Yichang,  P. R. China}
\email{changlin.xiang@ctgu.edu.cn}

\thanks{*Corresponding author: Chang-Lin Xiang}
\thanks{C.-Y. Guo is supported by the Young Scientist Program of the Ministry of Science and Technology of China (No.~2021YFA1002200), the National Natural Science Foundation of China (No.~12101362) and the Natural Science Foundation of Shandong Province (No.~ZR2022YQ01, ZR2021QA003). The corresponding author C.-L. Xiang is financially supported by the National Natural Science Foundation of China (No.~12271296 and 11701045).}

\begin{abstract}
In this survey paper, we give an overview of the conservation law approach in the study of geometric PDEs that models in particular polyharmonic maps.
\end{abstract}

\maketitle

{\small
\keywords {\noindent {\bf Keywords:} Elliptic system, Conservation law, Compensated compactness, Hodge decomposition, Gauge transform}
\smallskip
\newline
\subjclass{\noindent {\bf 2010 Mathematics Subject Classification:} 35J60, 53A10}}
\bigskip

\section{Motivation}
Interesting geometric partial differential equations are usually of critical or even supercritical nonlinearity in nature. The absence of possible applications of the maximum principle to solutions to non-linear elliptic systems reduces drastically the tools available for the regularity theory of weak solutions. In this aspect, the conservation law approach seems to be an effective strategy to address the regularity issues.

Consider the weakly harmonic map $u$ from the $n$-dimensional unit ball $B^n$ into the sphere $\mathbb{S}^{m}$ of $\R^{m+1}$. Chen \cite{C89} and Shatah \cite{Shatah-1988-CPAM} independently found a conservation law for $u$. Based on this conservation law, H\'elein \cite{H90} proved his celebrated regularity theorem for weakly harmonic maps from the two dimensional ball (or more generally surfaces) to sphere. Later, he also succeeded in extending the same result for general manifold targets by introducing the so-called moving frame technique; see \cite{Helein-2002}. The method of H\'elein is beautiful, but it does not apply for critical point of general second order conformally invariant elliptic Lagrangian with quadratic growth.

A major breakthrough towards the regularity issues of general second order conformally invariant elliptic variational problems was made by Rivi\`ere \cite{Riviere-2007}. He introduced the following second order linear system
\begin{equation}\label{eq:Riviere}
	-\Delta u=\Omega\cdot \nabla u
\end{equation}
and verified that \eqref{eq:Riviere} includes the Euler-Lagrange equation of general conformally invariant second order elliptic variational problems with quadratic growth in dimension two. Relying on (a variant of) the gauge theory of Uhlenbeck \cite{Uhlenbeck-1982}, Rivi\`ere succeeded in finding a conservation law for \eqref{eq:Riviere} and then regularity follows by standard potential theory. It should be noticed that this significant work not only gave a new proof of the regularity theorem of H\'elein, but also verifies affirmatively two long-standing regularity conjectures of Hildebrandt and Heinz on conformally invariant geometrical problems and the prescribed bounded mean curvature equations respectively; see \cite{Riviere-2007} for details.

The conservation law approach of \cite{Riviere-2007} was soon extended further by Rivi\`ere \cite{Riviere-2008} to Willmore surfaces, and by Lamm-Rivi\`ere \cite{Lamm-Riviere-2008} to fourth order linear elliptic system in dimension four, allowing him to give a new proof of the continuity (but not the stronger H\"older continuity) theorems of Chang et al. \cite{Chang-W-Y-1999} and Wang \cite{Wang-2004-CPAM} for biharmonic maps. Very recently, the conservation law approach was successfully extended by de Longueville and Gastel \cite{deLongueville-Gastel-2019} to general even order linear elliptic systems of Rivi\`ere type in the conformal dimension, providing a new proof of the continuity (but not the stronger H\"older continuity) result of polyharmonic maps \cite{Gastel-Scheven-2009CAG}; see also \cite{Horter-Lamm-2020} for a different construction of conservation law. Based on the above mentioned conservation law, Guo, Xiang and Zheng also established the H\"older continuity or indeed even an optimal $L^p$-regularity theory for fourth order system in \cite{Guo-Xiang-Zheng-2020-Lp} and for general even order system in \cite{Guo-Xiang-Zheng-2021-Lp}; thus yielding a complete recover of the regularity theorems of Chang et al. \cite{Chang-W-Y-1999}, Wang \cite{Wang-2004-CPAM} and Gastel-Scheven \cite{Gastel-Scheven-2009CAG}. Further more, in a recent work \cite{Guo-Xiang-Zheng-2022-CL}, Guo, Xiang and Zheng refined the conservation law of Lamm-Rivi\`ere \cite{Lamm-Riviere-2008} and \cite{deLongueville-Gastel-2019} so that they are fully equivalent to the equation on the whole domain of definition.

In the very recent work \cite{Guo-Xiang-2022-CL}, Guo and Xiang also partially extended the conservation law of Rivi\`ere \cite{Riviere-2007} for \eqref{eq:Riviere} to supercritical dimensions. In below, we shall give a detailed survey of the conservaiton law approach.

\section{The conservation law of Chen and Shatah}

Recall that a map $u\in W^{1,2}(B^n,\mathbb{S}^m)$ is weakly harmonic if it is a weak solution of 
\begin{equation}\label{eq:weakly harmonic into sphere}
	-\Delta u=|\nabla u|^2u.
\end{equation}
Because of the conformal invariance of the Dirichlet energy in dimension two, the equation \ref{eq:weakly harmonic into sphere} is also conformally invariant: if $u\colon B^2\to \mathbb{S}^m$ is a solution of \eqref{eq:weakly harmonic into sphere} and $\varphi\colon \R^2\to \R^2$ is a conformal map, then $u\circ \varphi$ is again a solution of \eqref{eq:weakly harmonic into sphere}. It is easily observed that the conformal dimension two is critical for \eqref{eq:weakly harmonic into sphere}: the right hand side of \eqref{eq:weakly harmonic into sphere} is merely in $L^1$ and so the classical $L^p$ regularity theory does not apply here. The following important conservation law was discovered, independently, by Chen \cite{C89} and Shatah \cite{Shatah-1988-CPAM}.

\begin{theorem}[Chen~\cite{C89}, Shatah~\cite{Shatah-1988-CPAM}]\label{thm:CL of Chen}
	A map $u\in W^{1,2}(B^2,\mathbb{S}^m)$ is weakly harmonic if and only if it satisfies the following conservation law
	\begin{equation}\label{eq:conservation law for sphere}
		\dive(u^i\nabla u^j-u^j\nabla u^i)=0\qquad \text{for all }i,j\in \{1,\cdots,m+1\}.
	\end{equation}
\end{theorem}

The conservation law of Chen and Shatah has a nice interpretation by means of Noether's theorem on conservation law: the symmetry (rotational invariance) of the target manifold $\mathbb{S}^{m}$ leads to \eqref{eq:conservation law for sphere}; see \cite{Helein-2002} for details.

We have the following celebrated regularity result due to H\'elein \cite{H90}.
\begin{theorem}[H\'elein]\label{thm:Helein sphere}
	If $u\in W^{1,2}(B^2,\mathbb{S}^m)$ is weakly harmonic, then it is smooth.
\end{theorem}
\begin{proof}
	By the conservation law of Chen and Shatah, we have
	\begin{equation*}
		\dive(u^i\nabla u^j-u^j\nabla u^i)=0\qquad \text{for all }i,j\in \{1,\cdots,n\}.
	\end{equation*}
	By the Hodge decomposition (or Poincar\'e's lemma), we know that there exists $B_{ij}\in W^{1,2}(B^2)$ such that $u^i\nabla u^j-u^j\nabla u^i=\nabla^{\perp}B_{ij}$, where $\nabla^{\perp}=(-\partial_y,\partial_x)$.
	
	Note that $|u|^2=1$ implies that $u^j\nabla u^j\footnote{Einstein's summation convention is used throughout this paper}=0$ and so the weakly harmonic map equation \eqref{eq:weakly harmonic into sphere} can be rewritten as 
	$$-\Delta u^i=(u^i\nabla u^j-u^j\nabla u^i)\cdot \nabla u^j=\nabla^\perp B_{ij}\cdot \nabla u^j.$$ 
	The product curl-grad on the right hand side has a Jacobian structure so that the classical result of Wente \cite{Went-1969} allows us to conclude the continuity of $u$.  
\end{proof}

\begin{remark}\label{rmk:on conservation law}
	i) Note that 
	$$\nabla^\perp B_{ij}\cdot \nabla u^j=-\nabla B_{ij}\cdot \nabla^{\perp}u^j=-\dive(B_{ij}\nabla^\perp u^j)$$
	and so we have the following conservation law
	$$-\dive(\nabla u^i-B_{ij}\nabla^{\perp}u^j)=0,$$
	or equivalently
	$$-\dive(\nabla u-B\nabla^{\perp}u)=0\quad \text{in }B^2.$$
	
	ii) Suppose $N\subset \R^m$ is a hypersurface and let $\nu\colon N\to T^\perp N$ be the smooth unit normal vector field. Given a weakly harmonic map $u\colon B^2\to N$, consider the composition function $w=\nu\circ u$. Using the fact that $\omega\cdot \partial_j u=0$, we may rewrite the harmonic map equation as
	\begin{equation}\label{eq:rewrite harmonic bad case}
		-\Delta u^i=w^i\partial_jw^k\partial_ju^k=\big(w^i\partial_jw^k-w^k\partial_jw^i \big)\partial_ju^k.
	\end{equation}
	Unlike the case $N=\mathbb{S}^m$, the vector fields
	$$
	v^{ik}:=\big(w^i\partial_jw^k-w^k\partial_jw^i  \big)_{1\leq j\leq 2}
	$$
	is no longer divergence free. Thus the preceeding proof fails to apply for general closed target manifold. 
	
	As we shall see, following H\'elein and Rivi\`ere, in spite of this apparent failture of the method, it is possible to once again produce the desired structure after a suitable transformation; in fact, as observed by Rivi\`ere \cite{Riviere-2007}, this is possible via gauge theory in the spirit of Uhlenbeck \cite{Uhlenbeck-1982}. For this, only the anti-symmetry $v^{ik}=-v^{ki}$ is needed.
\end{remark}

\section{The conservation law of Rivi\`ere }

In the revolutionary work \cite{Riviere-2007}, Rivi\`ere succeeded in finding a conservation law for the linear system \eqref{eq:Riviere}. More precisely, he proved the following result.

\begin{theorem}[Conservation law, \cite{Riviere-2007}]\label{thm:Riviere's conservation law}
	Suppose $\Omega\in L^2(B^2,so_m\otimes \Lambda^1\R^2)$. If there exist $A\in W^{1,2}\cap L^\infty(B^2,GL(m))$ and $B\in W^{1,2}(B^2,M(m))$ such that 
	\begin{equation}\label{eq:A B Omega}
		\nabla A+\nabla^\perp B=A\Omega,
	\end{equation}
	then $u$ solves \eqref{eq:Riviere} if and only if the following conservation law holds:
	\begin{equation}\label{eq:Riviere's CL}
		\dive(A\nabla u-B\nabla^\perp u)=0.
	\end{equation}
\end{theorem}
\begin{proof}
	This follows rather directly from computation:
	\[
	\begin{aligned}
		\dive(A\nabla u-&B\nabla^\perp u)=\nabla A\cdot \nabla u+A\Delta u-\nabla B\cdot \nabla^\perp u\\
		&=(\nabla A+\nabla^\perp B)\cdot \nabla u+A\Delta u\\
		&=A\Omega\cdot\nabla u+A\Delta u=A(\Omega\cdot \nabla u+\Delta u).
	\end{aligned}
	\]
\end{proof}

The main difficulty is thus to find $A$ and $B$ as in Theorem \ref{thm:Riviere's conservation law} that satisfies \eqref{eq:A B Omega}.

\begin{theorem}[Construction of conservation law, \cite{Riviere-2007}]\label{thm:construction of conservation law}
	There exists an $\epsilon_0=\epsilon_0(m)>0$ such that if $\Omega\in L^2(B^2,so_m\otimes \Lambda^1\R^2)$ satisfies 
	\[
	\|\Omega\|_{L^2(B^2)}\leq \epsilon_0,
	\]
	then there exist $A\in W^{1,2}\cap L^\infty(B^2,GL(m))$ and $B\in W^{1,2}(B^2,M(m))$ such that \eqref{eq:A B Omega} holds. Further more,  we have
	\begin{equation}\label{eq:A B controlled by Omega}
		\|\nabla A\|_{L^2(B^2)}^2+\|\nabla B\|_{L^2(B^2)}^2+\|\dist(A,SO_m)\|_{L^\infty(B^2)}\leq C(m)\|\Omega\|_{L^2(B^2)}^2.
	\end{equation}
\end{theorem}

The proof of Theorem \ref{thm:construction of conservation law} relies crucially on the following variant of the gauge theory of Uhlenbeck \cite{Uhlenbeck-1982}.

\begin{theorem}[Gauge transform, \cite{Riviere-2007}]\label{thm:gauge transform}
	Let $\Omega\in L^2(B^2,so_m\otimes \Lambda^1\R^2)$. Then there exist $\xi\in W_0^{1,2}(B^2,GL(m))$ and $P\in W^{1,2}(B^2,SO_m)$ such that
	\[
	P^{-1}\nabla P+P^{-1}\Omega P=\nabla^\perp \xi 
	\]
	and 
	\[
	\|\nabla^\perp \xi\|_{L^2(B^2)}^2+\|\nabla P\|_{L^2(B^2)}^2\leq C(m)\|\Omega\|_{L^2(B^2)}^2.
	\]
\end{theorem}
\begin{proof}[Sketch of the proof]
We follow the work of Schikorra \cite{S10}, who found a very elegant variational proof of this deep theorem. To be more precise, consider the variational problem
\[
\min_{Q\in W^{1,2}(B^2,SO_m)} E(Q)=\min_{Q\in W^{1,2}(B^2,SO_m)}\int_{B^2}|Q^{T}\nabla Q+Q^T\Omega Q|^2dx.
\]
One can show that there exists a minimizer $P\in W^{1,2}(B^2,SO_m)$ for the above functional, which satisfies the following Euler-Lagrange equation:
\[
\dive(P^T\nabla P+P^T\Omega P)=0.
\]
Green's formula then gives $\Omega_P\cdot \nu=0$, where $\Omega_P:=P^T\nabla P+P^T\Omega P$. Standard Hodge decomposition gives the existence of $\xi\in W^{1,2}_0(B^2,GL(m))$ such that
	\[
P^{-1}\nabla P+P^{-1}\Omega P=\nabla^\perp \xi. 
\]
The desired estimates follows easily from the minimizing property of $P$. For details, see \cite[Theorem 2.1]{S10}. 
\end{proof}

Using Theorem \ref{thm:gauge transform}, we now give a proof of Theorem \ref{thm:construction of conservation law}. The idea of the proof is due to Rivi\`ere \cite{Riviere-2007}, but with some minor technical improvement from \cite{Guo-Xiang-Zheng-2022-CL}.

\begin{proof}[Proof of Theorem \ref{thm:construction of conservation law}]
	By Theorem \ref{thm:gauge transform}, there exist $\xi\in W_0^{1,2}(B^2,GL(m))$ and $P\in W^{1,2}(B^2,SO_m)$ such that
	\[
	P^{-1}\nabla P+P^{-1}\Omega P=\nabla^\perp \xi 
	\]
	and 
	\[
	\|\nabla^\perp \xi\|_{L^2(B^2)}^2+\|\nabla P\|_{L^2(B^2)}^2+\|\nabla P^{-1}\|_{L^2(B^2)}^2\leq C(m)\|\Omega\|_{L^2(B^2)}^2.
	\]
	
	Observe that there exists a pair $(A,B)$ solves \eqref{eq:A B Omega} if and only if for $(\hat{A},B)$ with $\hat{A}=AP$, we have
	\begin{equation}\label{eq:hat A B P}
		\nabla \hat{A}-\hat{A}\nabla^{\perp}\xi=\nabla^{\perp}BP.
	\end{equation}
	It is thus sufficient to solve \eqref{eq:hat A B P} in $B^2$.
	
	We now use an extension argument from \cite{Guo-Xiang-Zheng-2022-CL} as follows: extend $\xi$ and $P$ from $B^2$ to $B_{2}^2$ such that $\xi=0$ and $P=I$ in $B_2^2\backslash \overline{B_{3/2}^2}$ (here we keep the notation $P$ and $\xi$ for the extended functions). Furthermore, we require the  norms of $P, \nabla P, \nabla \xi$ in $B_2^2$ is controlled by a constant multiple of the corresponding norms in $B^2$. Our strategy is to solve \eqref{eq:hat A B P} in the enlarged region $B_2^2$. 
	
	\textbf{Claim}. There exist $\bar{A}\in W^{1,2}(B_2^2,GL(m))$ and $B\in W^{1,2}_0(B_2^2,M(m))$ such that
	\begin{equation}\label{eq:auxillary equation}
		\begin{cases}
			\Delta \bar{A}=\nabla \bar{A}\cdot \nabla^{\perp}\xi+\nabla^{\perp}B\cdot \nabla P\\
			\Delta B=-\nabla^{\perp}\bar{A}\cdot \nabla P^{-1}-\dive\big((\bar{A}+I)\nabla\xi\cdot P^{-1} \big)\\
			\frac{\partial \bar{A}}{\partial \nu}=0 \text{ on}\quad \partial B^2_2 \quad \text{and}\quad \int_{B_2^2}\bar{A}=0\\
			B=0 \quad \text{on}\quad \partial B_2^2.
		\end{cases}
	\end{equation}
	Moreover, 
	\[
	\|\nabla \bar{A}\|_{L^2}+\|\bar{A}\|_{L^\infty}+\|\nabla B\|_{L^2}\leq C_m\|\Omega\|_{L^2(B^2)}.
	\]
	
	\textbf{Proof of Claim}. Set 
	\[
	X:=\{(a,b)\in W^{1,2}\cap L^\infty(B_2^2,GL(m))\times W^{1,2}(B_2^2,M(m)):\|(a,b)\|_X\leq 1 \},
	\]
	where $\|(a,b)\|_X:=\|\nabla a\|_{L^2}+\|\nabla b\|_{L^2}+\|a\|_{L^\infty}$.
	
	Standard elliptic regularity theory implies that for each $(a,b)\in X$, there exists a unique solution $(c,d)\in W^{1,2}\cap L^\infty\times W^{1,2}$ such that
	\[
	\begin{cases}
		\Delta c=\nabla a\cdot \nabla^{\perp}\xi+\nabla^{\perp}b\cdot \nabla P\\
		\Delta d=-\nabla^{\perp}a\cdot \nabla P^{-1}-\dive\big((a+I)\nabla\xi\cdot P^{-1} \big)\\
		\frac{\partial c}{\partial \nu}=0 \text{ on}\quad \partial B^2_2 \quad \text{and}\quad \int_{B_2^2}c=0\\
		d=0 \quad \text{on}\quad \partial B_2^2.
	\end{cases}
	\]
	Moreover, by Wente's lemma \cite{Went-1969}, for some $C=C_m>0$, we have
	\[
	\begin{aligned}
		\|\nabla c\|_{L^2}+\|c\|_{L^\infty}&\leq C\Big(\|\nabla a\|_{L^2}\|\nabla \xi\|_{L^2}+\|\nabla b\|_{L^2}\|\nabla P\|_{L^2} \Big)\\
		\|\nabla d\|_{L^2}&\leq C\Big(\|\nabla a\|_{L^2}\|\nabla P^{-1}\|_{L^2}+\|a\|_{L^\infty}\|\nabla \xi\|_{L^2}+\|\nabla \xi\|_{L^2} \Big)
	\end{aligned}
	\]
	This implies if $\|\Omega\|_{L^2(B^2)}$ is sufficiently small (less than $\frac{1}{2C}$), then
	\[
	\|(c,d)\|_X\leq C\|\Omega\|_{L^2}\cdot \|(a,b)\|_X+C\|\Omega\|_{L^2}\leq 2C\|\Omega\|_{L^2(B^2)}<1,
	\]
	which means $(c,d)\in X$. If we define 
	\[
	\begin{aligned}
		T\colon X&\to X\\
		(a,b)&\mapsto T(a,b):=(c,d),
	\end{aligned}
	\]
	then $T\colon X\to X$ is a contraction map. The fixed point theorem then gives a solution $(\bar{A},B)\in X$ that solves \eqref{eq:auxillary equation}.
	
	Set $\hat{A}:=A+I$. Then $(\hat{A},B)$ solves
	\begin{equation}\label{eq:aime equation}
		\begin{cases}
			\Delta \hat{A}=\nabla \hat{A}\cdot \nabla^{\perp}\xi+\nabla^{\perp}B\cdot \nabla P\\
			\Delta B=-\nabla^{\perp}\hat{A}\cdot \nabla P^{-1}-\dive\big(\hat{A}\nabla\xi\cdot P^{-1} \big)\\
			\frac{\partial \hat{A}}{\partial \nu}=0 \text{ on}\quad \partial B^2_2 \quad \text{and}\quad \int_{B_2^2}\hat{A}=4\pi I\\
			B=0 \quad \text{on}\quad \partial B_2^2.
		\end{cases}
	\end{equation}
	Moreover, we have
	\begin{equation}\label{eq:estimate on A B}
		\|\nabla\hat{A}\|_{L^2}+\|\hat{A}-I\|_{L^\infty}+\|\nabla B\|_{L^2}\leq C_m\|\Omega\|_{L^2(B^2)}.
	\end{equation}
	Observe that since $P=I$ and $\xi=0$  in $B_2^2\backslash \overline{B_{3/2}^2}$, we have $\bar{A}=0=B$ in $B_2^2\backslash \overline{B_{3/2}^2}$ and so $\hat{A}=I$ and $B=0$ in $B_2^2\backslash \overline{B_{3/2}^2}$.

	It remains to show \eqref{eq:hat A B P} holds for the pair $(\hat{A},B)$ in $B_2^2$ and thus in $B^2$ as well.
	
	Note first that by \eqref{eq:aime equation}, we have
	$$\dive(\nabla \hat{A}-\hat{A}\nabla^\perp \xi-\nabla^\perp BP)=0\quad \text{in }B_2^2.$$
	The Hodge decomposition implies that there exists $C\in W^{1,2}(B_2^2,M(m))$ such that
	\[
	\nabla \hat{A}-\hat{A}\nabla^\perp \xi-\nabla^\perp BP=\nabla^\perp C\quad \text{in }B_2^2.
	\]
	Our aim is to show $C\equiv 0$ in $B_2^2$. 
	
	By the second equation in \eqref{eq:aime equation}, we have $\dive(\nabla CP^{-1})=0$ and so Hodge decomposition implies that there exists $D\in W^{1,2}(B_2^2,M(m))$ such that 
	$$\nabla CP^{-1}=\nabla^\perp D\quad \text{in }B_2^2.$$
	Since $C=0$ in $B_2^2\backslash \overline{B_{3/2}^2}$, so is $D$. We may additionally assume $\int_{B_2^2}D=0$. Then we have by Wente's lemma again 
	\[
	\|\nabla D\|_{L^2}\leq C\|\nabla P^{-1}\|_{L^2}\|\nabla C\|_{L^2}\leq C\|\nabla D\|_{L^2}\|\Omega\|_{L^2(B^2)}\leq \frac{1}{2}\|\nabla D\|_{L^2},
	\]
	if $\|\Omega\|_{L^2}$ is sufficiently small. This implies that $\nabla D\equiv 0$ in $B_2^2$ and so is $\nabla C$. Since $C=0$ in $B_2^2\backslash \overline{B_{3/2}^2}$, we thus conclude $C\equiv 0$ in $B_2^2$ as desired. The estimate \eqref{eq:A B controlled by Omega} follows directly from \eqref{eq:estimate on A B}.
\end{proof}

\begin{remark}\label{rmk:on construction of conservation law}
	It follows from the proof of Theorem \ref{thm:construction of conservation law} that there are indeed infinitely many choice of $A$ and $B$ as required, as there are infinitely many ways to extend $P$ and $\xi$ with the desired properties.
\end{remark}

\section{The conservation law for higher order system of Rivi\`ere type}

In an interesting recent work \cite{deLongueville-Gastel-2019}, de Longueville and Gastel introduced the following even order linear elliptic system of Rivi\`ere type
\begin{equation}\label{eq:Longue-Gastel system}
	\Delta^{m}u=\sum_{l=0}^{m-1}\Delta^{l}\left\langle V_{l},du\right\rangle +\sum_{l=0}^{m-2}\Delta^{l}\delta\left(w_{l}du\right) \qquad \text{ in }B^{2m}.
\end{equation}
System \eqref{eq:Longue-Gastel system} includes (both extrinsic and intrinsic) $m$-polyharmonic mappings. It reduces to the Lamm-Rivi\`ere system \cite{Lamm-Riviere-2008} when $m=2$, and to \eqref{eq:Riviere} when $m=1$. The coefficient functions are assumed to satisfy
\begin{equation}\label{eq:coefficient w V}
	\begin{aligned}
		&w_{k} \in W^{2 k+2-m, 2}\left(B^{2 m}, \mathbb{R}^{n \times n}\right) \quad \text { for } k \in\{0, \ldots, m-2\} \\
		&V_{k} \in W^{2 k+1-m, 2}\left(B^{2 m}, \mathbb{R}^{n \times n} \otimes \wedge^{1} \mathbb{R}^{2 m}\right) \quad \text { for } k \in\{0, \ldots, m-1\}.
	\end{aligned}
\end{equation}
Moreover,  the first order potential $V_0$ has the decomposition $V_{0}=d \eta+F$ with
\begin{equation}\label{eq:coefficient eta F}
	\eta \in W^{2-m, 2}\left(B^{2 m}, s o(n)\right), \quad F \in W^{2-m, \frac{2 m}{m+1}, 1}\left(B^{2 m}, \mathbb{R}^{n \times n} \otimes \wedge^{1} \mathbb{R}^{2 m}\right).
\end{equation}

To formulate the conservation law, we set
\begin{equation}\label{eq:theta for small coefficient}
	\begin{aligned}
		\theta_{D}:=\sum_{k=0}^{m-2}&\|w_k\|_{W^{2k+2-m,2}(D)}+\sum_{k=1}^{m-1}\|V_k\|_{W^{2k+1-m,2}(D)}\\
		&+\|\eta\|_{W^{2-m,2}(D)}+\|F\|_{W^{2-m,\frac{2m}{m+1},1}(D)}
	\end{aligned}
\end{equation}
for $D\subset \R^{2m}$.

We have the following conservation law for \eqref{eq:Longue-Gastel system}.
\begin{theorem}[Conservation law, \cite{Guo-Xiang-Zheng-2022-CL}]\label{thm:general even order} 
	There exist  constants $\epsilon_{m},C_m>0$ such that under the smallness assumption $\theta_{B^{2m}}<\epsilon_m$, there exist $A\in W^{m,2}\cap L^\infty(B^{2m},Gl(n))$ and $B\in W^{2-m,2}(B^{2m},\R^{n\times n}\otimes \wedge^2\R^{2m})$ satisfying 
	\begin{equation}\label{eq:for CL general order}
		\Delta^{m-1}dA+\sum_{k=0}^{m-1}(\Delta^k A)V_k-\sum_{k=0}^{m-2}(\Delta^k dA)w_k=\delta B\qquad \text{in }B^{2m}.
	\end{equation} 
	Moreover, 
	\[
	\|A\|_{W^{m,2}(B^{2m})}+\|dist(A,SO(m))\|_{L^{\infty}(B^{2m})}+\|B\|_{W^{2-m,2}(B^{2m})}\leq C_{m}\theta_{B^{2m}}.
	\]
	Consequently, $u$ solves \eqref{eq:Longue-Gastel system} if and only if it satisfies the conservation law
	\begin{equation}\label{eq:conservation law of D-G}
		\begin{aligned}
			0&=\delta\Big[\sum_{l=0}^{m-1}\left(\Delta^{l} A\right) \Delta^{m-l-1} d u-\sum_{l=0}^{m-2}\left(d \Delta^{l} A\right) \Delta^{m-l-1} u \\ &\qquad -\sum_{k=0}^{m-1} \sum_{l=0}^{k-1}\left(\Delta^{l} A\right) \Delta^{k-l-1} d\left\langle V_{k}, d u\right\rangle+\sum_{k=0}^{m-1} \sum_{l=0}^{k-1}\left(d \Delta^{l} A\right) \Delta^{k-l-1}\left\langle V_{k}, d u\right\rangle \\ &\qquad -\sum_{k=0}^{m-2} \sum_{l=0}^{k-2}\left(\Delta^{l} A\right) d \Delta^{k-l-1} \delta\left(w_{k} d u\right)+\sum_{k=0}^{m-2} \sum_{l=0}^{k-2}\left(d \Delta^{l} A\right) \Delta^{k-l-1} \delta\left(w_{k} d u\right) \\ &\qquad -\langle B, d u\rangle\Big]
		\end{aligned}
	\end{equation} 
	in $B^{2m}$.  
\end{theorem}

Theorem \ref{thm:general even order} was proved in \cite{Lamm-Riviere-2008} for the case $m=2$ and in \cite{deLongueville-Gastel-2019} for general $m$, but with the conservation law \eqref{eq:conservation law of D-G} holding only on $B_{1/2}^{2m}$. An extra extension argument as was done in the proof of Theorem \ref{thm:construction of conservation law} would lead to the current form; see \cite{Guo-Xiang-Zheng-2022-CL} for details.

\section{Open problems}
%

In the final section, we shall discuss some natural open problems related to the conservation law. 

Rivi\`ere's theory about \eqref{eq:Riviere} is very successful in the critical dimension $n=2$, but it is much less understood in higher dimensions. In \cite[Page 7]{Riviere-2007}, Rivi\`ere pointed out that
\begin{problem}\label{prob:1}
	Can we establish the conservation law for \eqref{eq:Riviere} in higher dimensions $n\geq 3$? 
\end{problem}

In \cite{Riviere-2007}, Rivi\`ere proposed to look for conservation laws in the natural Morrey space $M^{2,n-2}$, which has direct application to regularity of stationary harmonic maps in higher dimensions. Later, Rivi\`ere and Struwe \cite{Riviere-Struve-2008} counstructed an example showing that Wente's lemma fails with coefficients merely in $M^{2,n-2}$ and thus it is not easy to find a conservation law in such spaces. In the very recent work \cite{Guo-Xiang-2022-CL}, we succeeded in finding the conservation law in the smaller Lorentz space $L^{n,2}\subsetneq M^{2,n-2}$. This is, however, far from a satisfied theory in higher dimensions.

The second natural problem is
\begin{problem}\label{prob:2}
	What is the essence of Rivi\`ere's conservation law?  
\end{problem}

Unlike Noether's conservation law, there is no variational characterization of Rivi\`ere's conservation law. The existence of $A$ and $B$ appearing in \eqref{eq:A B Omega} are obtained through a fixed point type argument and thus are not explicit. Furthermore, as was pointed out in Remark \ref{rmk:on construction of conservation law}, there are in fact infinitely many $A$ and $B$ such that the conservation law \eqref{eq:Riviere's CL} holds. Thus it is natural to find a good explanation of this conservation law from a scientific point of view. 

Another fundamental problem in Rivi\`ere's theory or the theory of harmonic maps is whether one can establish the global regularity in higher dimensions. The following well-known conjecture was asked by Rivi\`ere in \cite[Page 9 Conjecture]{Riviere-2007}:
\begin{problem}\label{prob:3}
	Conjecture that for every $k\leq m$, for every $n\in \mathbb{N}$, for every $k$-dimensional closed submanifold $N$ of $\mathbb{R}^m$, and for every $C>0$, there exists $\delta=\delta(C,n,N)>0$ such that if $u$ is a $W^{1,2}$ weakly harmonic map from $B_2^n(0)$ into $N$ satisfying 
	\[
	\int_{B_2^n(0)}|\nabla u|^2dx\leq C,
	\]
	then 
	\[
	\int_{B_1^n(0)}|\nabla^2 u|dx\leq \delta. 
	\]
\end{problem}

Similar type of result has been obtained by Lin \cite{L99}, Naber and Valtorta \cite{NV17}, for stationary harmonic maps, but with extra topological restrictions on the target manifold $N$. Without any further requirement on the target manifold, the result was only known for minimizing harmonic maps. The usual small $\epsilon$-energy regularity result implies that the conjecture holds if $C$ is sufficiently small (quantitatively). For large $C$, there is no theory available. Note however that a positive answer to Problem \ref{prob:3} would imply the corresponding energy identity for stationary harmonic maps into genearl closed manifolds via the techniques of Lin-Rivi\`ere \cite{LR02}. The classical applicaiton of conservation law in regularity issues includes the $L^p$-estimates for inhomogeneous system of Rivi\`ere type; see \cite{Sharp-Topping-2013-TAMS,Guo-Xiang-Zheng-2020-Lp,Guo-Xiang-Zheng-2021-Lp}. Towards a possible solution of Problem \ref{prob:3}, it is suggestive to first establish a conservation law for harmonic maps in higher dimensions.

\medskip
\textbf{Acknowledgements.} We would like to thank Prof.~Armin Schikorra, who kindly pointed out to us the beautiful work \cite{S10}, which makes the proof of Theorem \ref{thm:construction of conservation law} essentially self-contained.

\end{document}